\documentclass{gtart_h}  

\def\ifplaintex{\expandafter\ifx\csname documentclass\endcsname\relax}

\ifplaintex 
\else
\fi

\expandafter\ifx\csname beginpicture\endcsname\relax
\expandafter\ifx\csname documentclass\endcsname\relax
\input pictex \else
\input prepictex \input pictex \input postpictex \fi\fi

\def\gt{{\mathsurround=0pt\it $\cal G\mskip-2mu$eometry \&\ 
$\cal T\!\!$opology}}        

\def\gtp{{\mathsurround=0pt\it $\cal G\mskip-2mu$eometry \&\ 
$\cal T\!\!$opology $\cal P\!$ublications}}  

\def\lognumber#1{\def\thelognumber{#1}}
\def\volumenumber#1{\def\thevolumenumber{#1}}
\def\papernumber#1{\def\thepapernumber{#1}}
\def\volumeyear#1{\def\thevolumeyear{#1}}

\def\pagenumbers#1#2{\def\startpage{#1}\def\finishpage{#2}}
\def\published#1{\def\publishdate{#1}}
\def\proposed#1{\def\theproposer{#1}}
\def\seconded#1{\def\theseconders{#1}}
\def\received#1{\def\receiveddate{#1}}

\def\accepted#1{\def\accepteddate{#1}}
\def\asciititle#1{\def\theasciititle{#1}}
\def\covertitle#1{\def\thecovertitle{#1}}

\long\def\asciiabstract#1{\long\def\theasciiabstract{#1}}
\def\asciikeywords#1{\def\theasciikeywords{#1}}

\let\\\par\let\thelognumber\relax
\let\thevolumenumber\relax\let\thepapernumber\relax
\let\thevolumeyear\relax\let\thesamplenumber\relax\let\startpage\relax
\let\finishpage\relax\let\publishdate\relax\let\receiveddate\relax
\let\reviseddate\relax\let\accepteddate\relax\let\theasciititle\relax
\let\thecovertitle\relax\let\theasciiauthors\relax
\let\theasciiabstract\relax\let\theasciikeywords\relax
\let\theasciiemail\relax\let\theshortauthors\relax\let\theshorttitle\relax

\long\def\maketitlep{   

\count0=\startpage

\gt\hfill      
\beginpicture
\setcoordinatesystem units <0.33truein, 0.33truein> point at 2.2 0.9
\setplotsymbol ({$\cal G$})
\plotsymbolspacing=9truept
\circulararc 315 degrees from 0 1 center at 0 0
\setplotsymbol ({$\cal T$})
\circulararc 315 degrees from 1 -1 center at 1 0
\endpicture
\break
{\small\ifx\thesamplenumber\relax 
Volume \else Sample
\fi\thevolumenumber\ (\thevolumeyear)
\startpage--\finishpage\nl
Published: \publishdate}
\vglue 0.5truein plus 0.4fil minus 0.1truein

{\parskip=0pt\leftskip 0pt plus 1fil\def\\{\par\smallskip}{\ifplaintex\large
\else\Large\fi\bf\thetitle}\par\medskip}   

\vglue 0pt plus 0.1fil 

{\parskip=0pt\leftskip 0pt plus 1fil\def\\{\par}{\sc\theauthors}
\par\medskip}

\vglue 0pt plus 0.1fil 

{\small\parskip=0pt\let\newline\\
{\leftskip 0pt plus 1fil\def\\{\par}{\sl\theaddress}\par}
\expandafter\ifx\theemail\relax    
\relax\else\vglue 5pt plus 0.02fil minus 2pt\def\\{\stdspace{\rm 
and}\stdspace} 
\cl{Email:\stdspace\tt\theemail}\fi
\ifx\theurl\relax                  
\relax\else\vglue 5pt plus 0.02fil minus 2pt\def\\{\stdspace{\rm 
and}\stdspace}
\cl{URL:\stdspace\tt\theurl}\fi\par}

\vglue 7pt plus 0.3fil minus 3pt

{\bf Abstract}
\vglue 5pt plus 0.1fil minus 2pt

\theabstract

\vglue 7pt plus 0.3fil minus 3pt

{\bf AMS Classification numbers}\quad Primary:\quad \theprimaryclass

Secondary:\quad \thesecondaryclass

\vglue 5pt plus 0.3fil minus 2pt

{\bf Keywords:}\quad \thekeywords

\vglue 10pt plus 0.5fil minus 5pt

{\small  Proposed: \theproposer\hfill Received: \receiveddate\nl
Seconded: \theseconders\hfill 
\ifx\reviseddate\relax                         
Accepted: \accepteddate                        
\else
Revised: \reviseddate                          
\fi}
\eject
}       

\let\maketitlepage\maketitlep
\let\maketitle\maketitlepage

\font\phead=cmsl9 scaled 950
\font=cmsl9 scaled 1050
\font\pnum=cmbx10 scaled 913
\font\lnum=cmbx10 
\font\pfoot=cmsl9 scaled 950
\font=cmsl9 scaled 1050
\ifplaintex
\headline{\vbox to 0pt{\vskip -4.5mm\line{\small\phead\ifnum
\count0=\startpage ISSN 1364-0380 (on line)
1465-3060 (printed) \hfill {\pnum\folio}\else\ifodd\count0\def\\{ }%
\ifx\theshorttitle\relax\thetitle\else\theshorttitle\fi\hfill{\pnum\folio}
\else\def\\{ and }{\pnum\folio}\hfill\ifx\theshortauthors\relax\theauthors
\else\theshortauthors\fi\fi\fi}\vss}}
\footline{\vbox to 0pt{\vglue 0mm\line{\small\pfoot\ifnum\count0=\startpage
\copyright\ \gtp\hfill\else
\gt, Volume \thevolumenumber\ (\thevolumeyear)\hfill\fi}\vss
}}
\else
\makeatletter
\def\@oddhead{{\small\lhead\ifnum\count0=\startpage ISSN 1364-0380 (on line)
1465-3060 (printed) \hfill {\lnum\number\count0}\else\ifodd\count0
\def\\{ }\ifx\theshorttitle\relax \thetitle \else\theshorttitle\fi\hfill
{\lnum\number\count0}\else\def\\{ and }{\lnum\number\count0}
\hfill\ifx\theshortauthors\relax 
\theauthors\else\theshortauthors\fi\fi\fi}}\def\@evenhead{\@oddhead}
\def\@oddfoot{\small\lfoot\ifnum\count0=\startpage\copyright\ \gtp\hfill\else
\gt, Volume \thevolumenumber\ (\thevolumeyear)\hfill\fi}
\def\@evenfoot{\@oddfoot}
\makeatother
\fi

\newwrite\gtoutfile
\long\gdef\makeheadfile{  
{\def\\{, }\def\s{ }
\immediate\openout\gtoutfile head.xxx
\immediate\write\gtoutfile{Proxy-for: \ifx\theasciiauthors\relax
\theauthors\else\theasciiauthors\fi\s<\ifx\theasciiemail\relax\theemail\else\theasciiemail\fi>}
\immediate\write\gtoutfile{\noexpand\\}
\immediate\write\gtoutfile{Authors: \ifx\theasciiauthors\relax
\theauthors\else\theasciiauthors\fi}
{\def\\{ }\immediate\write\gtoutfile{Title: \ifx\theasciititle\relax
\thetitle\else\theasciititle\fi}}
\immediate\write\gtoutfile{Subj-class: GT or SG or MG etc}
\immediate\write\gtoutfile{MSC-class: \theprimaryclass\ifx\thesecondaryclass\relax\else, \thesecondaryclass\fi}
\immediate\write\gtoutfile{Journal-ref: Geom. Topol. \thevolumenumber
(\thevolumeyear) \startpage-\finishpage}
\immediate\write\gtoutfile{Comments: Published by Geometry and Topology at}
\immediate\write\gtoutfile{\s\s http://www.maths.warwick.ac.uk/gt/GTVol\thevolumenumber/paper\thepapernumber.abs.html}
\immediate\write\gtoutfile{\noexpand\\}
\immediate\write\gtoutfile{}
\ifx\theasciiabstract\relax
\immediate\write\gtoutfile{\theabstract}\else
\immediate\write\gtoutfile{\theasciiabstract}\fi
\immediate\write\gtoutfile{}
\immediate\write\gtoutfile{\noexpand\\}
\immediate\write\gtoutfile{}
\immediate\closeout\gtoutfile}}  

\def\maketitlepage{\maketitlep\makeheadfile}
\let\maketitle\maketitlepage

\lognumber{425}
\volumenumber{8}\papernumber{19}\volumeyear{2004}
\pagenumbers{735}{742}
\received{20 Febrary 2004}
\accepted{29 April 2004}
\published{17 May 2004}
\proposed{Peter Ozsv\'ath}
\seconded{Joan Birman, Ronald Fintushel}

\usepackage{amssymb, amsmath, graphicx}

\newtheorem{theorem}{Theorem}

\newtheorem{corollary}[theorem]{Corollary}
 
\theoremstyle{definition}

\newtheorem{example}[theorem]{Example}

\newcommand{\fig}[2] { \includegraphics[scale=#1]{#2}  }
\def\zz{{\bf Z}}
\def\calc{\mathcal{C}}

\def\cala{\mathcal{A}}
\def\calp{\mathcal{P}}

\begin{document}

\title {Computations of the Ozsv\'ath--Szab\'o\\knot concordance invariant}
\covertitle{Computations of the Ozsv\noexpand\'ath--Szab\noexpand\'o\\knot concordance invariant}
\asciititle{Computations of the Ozsvath-Szabo knot concordance invariant}
\author{Charles Livingston}
\address{Department of Mathematics, Indiana University\\Bloomington, IN  47405,
USA}
\email{livingst@indiana.edu}

\begin{abstract}
Ozsv\'ath and Szab\'o have defined a  knot concordance invariant $\tau$ that
bounds the 4--ball genus of a knot.  Here we discuss shortcuts to its computation.  We include examples
of Alexander polynomial one knots   for which the invariant is nontrivial, including all iterated
untwisted positive doubles of  knots with nonnegative Thurston--Bennequin number, such as the
trefoil, and explicit computations for several 10 crossing knots. We also note that a new proof of the
Slice--Bennequin Inequality quickly follows from these techniques.
\end{abstract}
\asciiabstract{%
Ozsvath and Szabo have defined a knot concordance invariant tau
that bounds the 4-ball genus of a knot.  Here we discuss shortcuts to
its computation.  We include examples of Alexander polynomial one
knots for which the invariant is nontrivial, including all iterated
untwisted positive doubles of knots with nonnegative
Thurston-Bennequin number, such as the trefoil, and explicit
computations for several 10 crossing knots. We also note that a new
proof of the Slice-Bennequin Inequality quickly follows from these
techniques.}

\primaryclass{57M27}
\secondaryclass{57M25, 57Q60}
\keywords{Concordance, knot genus, Slice--Bennequin Inequality}
\asciikeywords{Concordance, knot genus, Slice-Bennequin Inequality}

\maketitle

Using their theory of knot Floer homology, Ozsv\'ath and Szab\'o \cite{os}  defined an   invariant
$\tau$ of knots in $S^3$ and showed that it induces a homomorphism
$\tau\co \calc \to
 \zz$, where $\calc$ is the concordance group of smooth knots in $S^3$. Computations of $\tau$ for
particular knots, and more generally the application of knot Floer homology to bound the 4--ball genus
of knots (eg \cite{ow, pl, ra1, ra2}),   depend upon a  detailed understanding of its definition.  Here we show
that the most basic properties of $\tau$ developed in
\cite{os} are sufficient to yield its quick evaluation  for a number of interesting examples including
some pretzel knots of Alexander polynomial one, iterated untwisted doubles of knots with nonnegative
Thurston--Bennequin number and some interesting 10 crossing knots. 

Although we do not use the deeper theoretical work of Rudolph (eg~\cite{r1}) here, in   ways our
approach parallels  his extension of the results of Kronheimer--Mrowka~\cite{km} on torus knots to more
general knots and his proof of the Slice--Bennequin Inequality.

Three essential properties of~$\tau$ are stated  in the following theorem.

\begin{theorem}\label{thmos} There exists an integer valued knot invariant $\tau$
  satisfying:
\begin{enumerate}

\item $\tau(K \# J) = \tau(K) + \tau(J)$ and $\tau(-K) = - \tau(K)$ for all knots $K$ and $J$.
\item  The value of
$\tau$ is bounded by the smooth 4--ball genus, $ \tau(K)  \le g_4(K)$.

\item For  the
$(p,q)$--torus knot with $p, q > 0$,  $T_{p,q}$,  $\tau$ equals the 3--sphere genus,
$g_3(T_{p,q})$. Specifically, $\tau(T_{p,q}) =
  (p-1)(q-1)/2 $.
\end{enumerate}
\end{theorem}

 An immediate consequence, as described in~\cite{os}, is:

\begin{corollary} $\tau$ induces a homomorphism $\tau\co \calc \to \zz$ and $ |\tau(K)|  \le g_4(K)$.
\end{corollary}

\begin{proof}That $ |\tau(K)|  \le g_4(K)$ follows from $\tau(-K) = - \tau(K)$, $g_4(K) = g_4(-K)$ and  $
\tau(K)  \le g_4(K)$.  Next, if
$K$ is concordant to $J$, then $K \# -J$ is slice, and hence of 4--genus 0.  Thus, $\tau(K) + \tau(-J) =
0$, $\tau(K) = \tau(J)$, and so
  $\tau $ is a concordance invariant.
\end{proof}

 The following   appears in~\cite{os} as a corollary of the general relationship between
$\tau(K)$ and the genus of surfaces bounded by $K$ in negative definite 4--manifolds. Here we note that
it follows immediately from Theorem~\ref{thmos}.

\begin{corollary} If $K_+$ and $K_-$ differ by a single crossing change, from positive to negative, then
$ 0 \le\tau(K_+) -\tau(K_-)  
\le 1$.
\end{corollary}

 \begin{proof}The crossing change provides a genus 1 cobordism from $K_+$ to
 $K_-$.  Thus, $g_4(K_+ \# {-}K_-)\le 1$ and so $|\tau(K_+) - \tau( K_-) | \le
 1$.  A negative crossing change converts $-T_{2,3}$ into the unknot, so $(K_+ \# -T_{2,3}) \# -K_-$
 bounds a disk in $B^4$ with two double points  of opposite signs.  Tubing these double points together
 shows that $g_4((K_+ \# -T_{2,3}) \# -K_-) \le 1$.  Thus, $|\tau(K_+) - 1 - \tau(K_-)|$\break$  \le 1$. 
 Combining the two inequalities gives the desired result.
 \end{proof}

\section{Subsurfaces of Torus Knot Fibers}

  For a surface $F$ we let $g(F)$ denote the genus of $F$.  Recall that for any torus knot $T_{p,q}$ the
complement is fibered over $S^1$ and  the fiber $F$ realizes the 3--genus of $T_{p,q}$, $ (p-1)(q-1)/2$.

 \begin{theorem}\label{mainthm} Suppose that a knot $K$ is embedded in the interior of a fiber surface
$F$ of a torus knot $  T = T_{p,q}$ with
$pq >0$ and that
$K$ is null homologous on
$F$, bounding a surface $G \subset F$.  Then  $\tau(K) = g_4(K) = g_3(K) = g(G)$. \end{theorem}

\begin{proof}  A Morse function $h \co F- \mbox{int}(G) \to  [0,1]$  taking value 0 on
$K$ and
$1$ on $T$ gives the cobordism $(\mbox{id} \times h) \co  F- \mbox{int}(G)
\to  S^3 \times [0,1]$ from $T_{p,q}$ to $K$ of genus $g(F) - g(G)$.  Hence,
$T
\# {-}K$ bounds a surface of genus $g(F) - g(G)$ in $B^4$ and
$g_4(T
\# {-} K) \le g(F) - g(G)$.   Thus $\tau(T) - \tau(K) \le g(F) - g(G)$. By Theorem~\ref{thmos},
$\tau(T) = g(F)$ and hence
$g(G)
\le
\tau(K)$. We then have the string of inequalities
$$\tau(K) \le g_4(K)
\le g_3(K)
\le g(G) \le \tau(K)$$ and these yield the desired result.
\end{proof}

In the language of Rudolph (eg~\cite{r1}), such surfaces $G$ on fibers of torus knots are called {\it
quasipositive surfaces}. One quick consequence of Rudolph's work is the following.

\begin{corollary} The  untwisted positive double of the trefoil and the pretzel knot
$P(3,-5,-7)$ (both of Alexander polynomial one) have $\tau = 1$.  
\end{corollary}
\begin{proof} Rudolph~\cite{r1} has drawn an explicit illustration of  
$P(3,-5,-7)$ on the fiber surface for the torus link $T_{5,5}$ and that illustration applies as well
for $T_{5,6}$.   Rudolph also indicates how a similar illustration can be drawn for the double of the
trefoil.  (Our sign convention for pretzel knots here is the opposite of that in~\cite{r1} and is consistent with Rudolph's more general work on pretzel knots in~\cite{r3}.  The positive double is the double formed from two parallel unlinked copies of the knot by adding a clasp with two crossing points, both of which have positive sign.)
\end{proof}

\begin{corollary} The subgroup $\calp \subset \calc$ generated by knots of Alexander polynomial one
contains a summand isomorphic to $\zz$.  The knot $P( 3,-5,-7)$ represents a generator of such a
summand  and in particular is not divisible: $P(3,-5,-7) \ne aK \in \calc$ for any $a \ne \pm 1$. 
\end{corollary}

\begin{proof}The argument is the same as in~\cite{l} where $\cala$ instead of $\calp$ is considered. 
Since $\tau$ maps $\calp$ onto the free abelian group $\zz$ the map splits.  Since $\tau(P(3,-5,-7)) =
1$, $P(3,-5,-7)$ cannot be divisible.
\end{proof}

Rudolph's results along these lines have further applications.  In particular, his work on pretzel
knots~\cite{r2} implies  that the pretzel knot 
$P(t_1,
\ldots , t_k)$ with all $t_i$ and $k$ odd (and its standard Seifert surface) embeds on the fiber of a
torus knot with $pq > 0$ if and only if $t_i + t_j <0 $ for all $1 \le i<j \le k$ and thus for such pretzel knots
$\tau = (k-1)/2$.

\begin{corollary}\label{posbraid} If $\hat{\beta}$ is the closure of a positive braid of
$n$ strands and word length $k$, then $\tau(\hat{\beta}) = \frac{k-n+1}{2} = g_4(\hat{\beta})$.
\end{corollary}

\begin{proof} The torus knot $T_{n,q}$ can be drawn as an $n$--stranded positive braid.  Its fibered Seifert surface is formed from $n$ parallel disks joined by $q(n-1)$ twisted bands, one for each crossing point.  Similarly, the  Seifert surface $G$ for $\hat{\beta}$ is built from $n$ parallel disks by joining them
with $k$ twisted bands. The resulting surface has Euler characteristic $n - k$, and hence genus
$\frac{k-n+1}{2}$.  By adding more bands to this surface, one can construct the fiber Seifert surface
$F$ of the torus knot $T_{n,q}$ for some large
$q$. Thus, removing a  small open tubular neighborhood, on $G$, of the boundary of $G$ yields a surface
homeomorphic to $G$ with boundary  isotopic to $K$ in the interior of
$F$.  The proof is now completed by applying Theorem~\ref{mainthm}.
\end{proof}

\section{Examples}

In~\cite{ka} there is a table listing the 4--genera of prime knots with 10 or fewer crossings, as then
known.  Since the appearance of that table, most of the unknown values have been determined. 
(See~\cite{sh} for an updated table, where it appears that $10_{51}$ remains as the only unknown
case.)  However, the unknown cases of~\cite{ka} continue to provide interesting test cases for new
techniques, since classical methods could not resolve them.  The few examples presented here were
chosen because of their appearance in~\cite{ow, os}, where direct application of the theory therein
developed was used.  Explicit calculations were not included in~\cite{os}, in that they were lengthy and
demanded the use of Mathematica.

\begin{figure}[ht!]
\centerline{  \fig{.25}{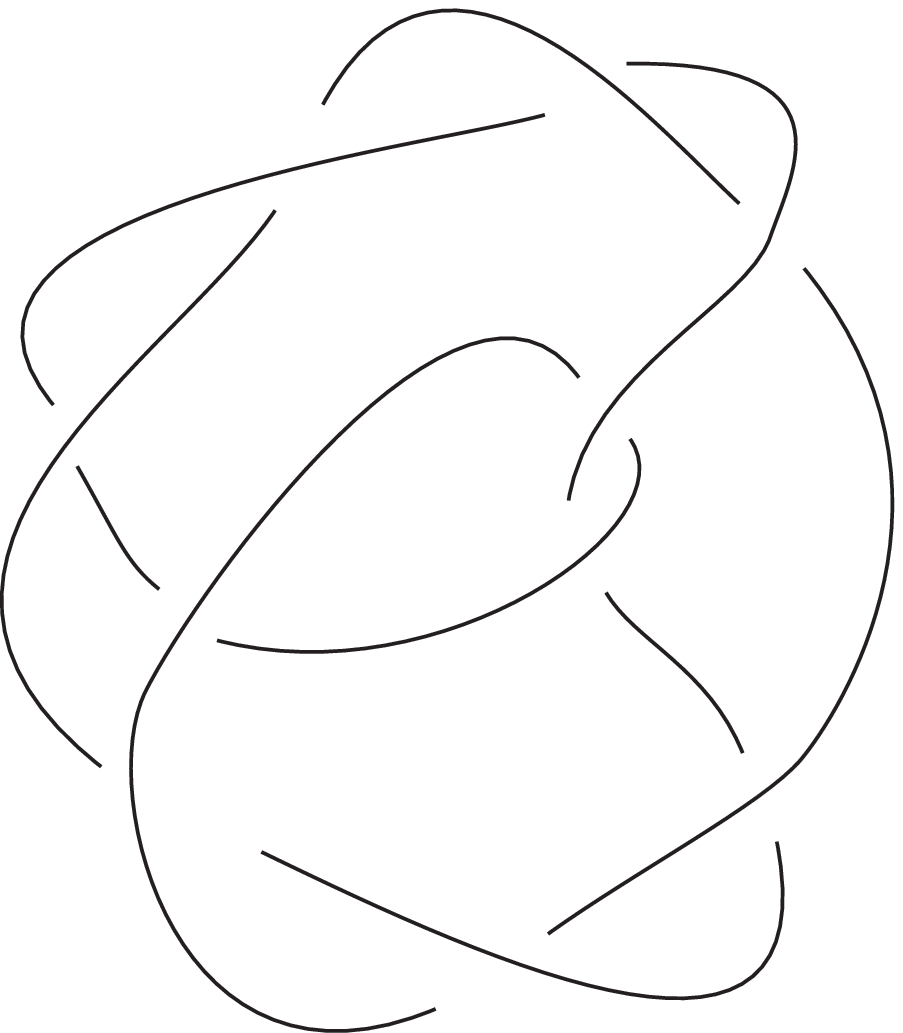}\hskip1in  \fig{.20}{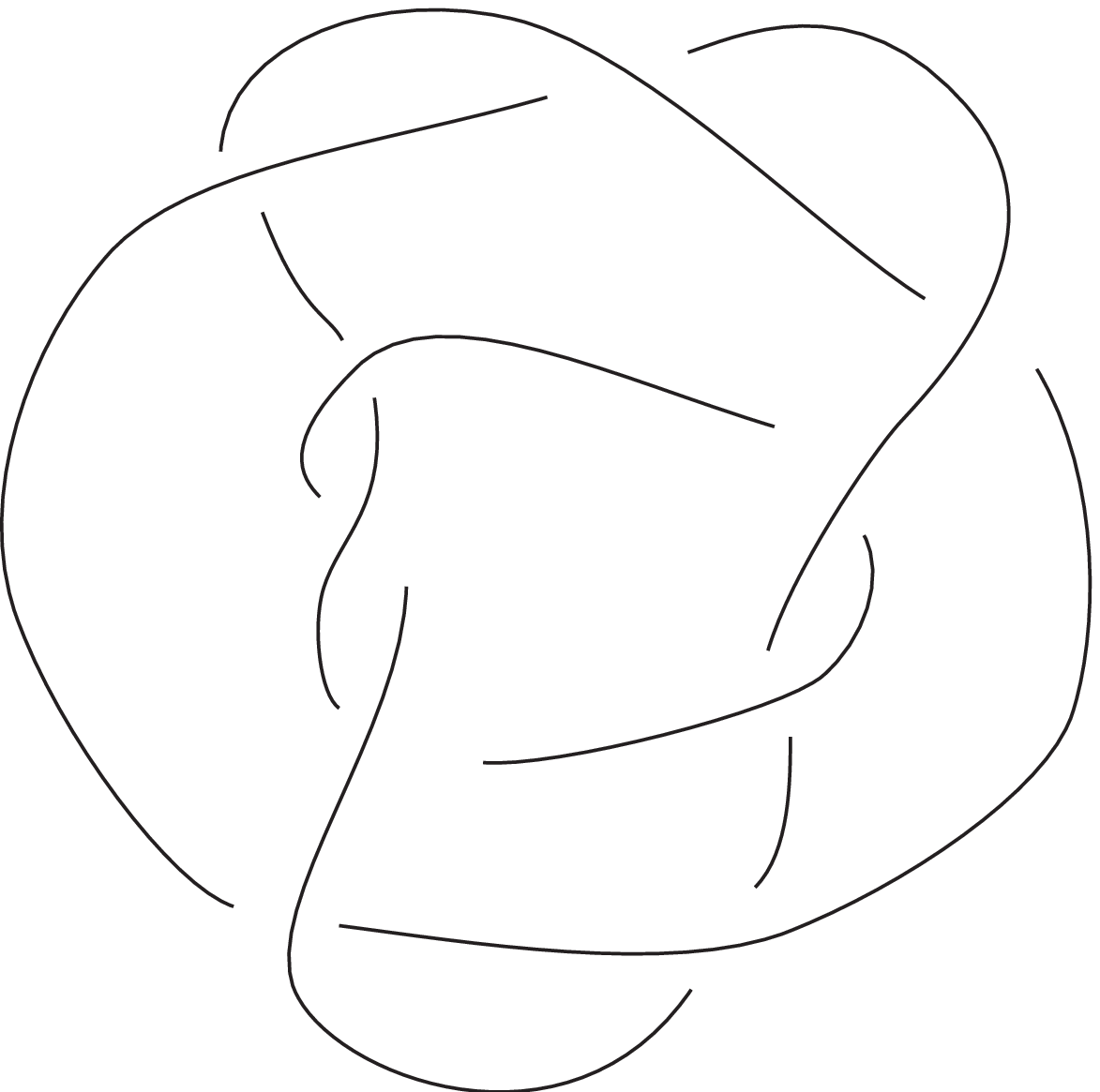}  }
 \caption{$10_{139}$ and  --$10_{152}$} \label{fig10_139}
\end{figure}

\begin{example} Using the notation of Rolfsen~\cite{rolf}, the knots
$10_{139}$ and  
${-}10_{152}$ as shown in Figure~\ref{fig10_139} each are closures of positive braids with 3 strands and
word length 10.  Thus, by  Corollary~\ref{posbraid}, the value of
$\tau$ and $g_4$ for each of these is 4.  (The 4--ball genus of each of these was first computed in~\cite{kawam}.) \end{example}

\begin{figure}[ht!]
\centerline{  \fig{.25}{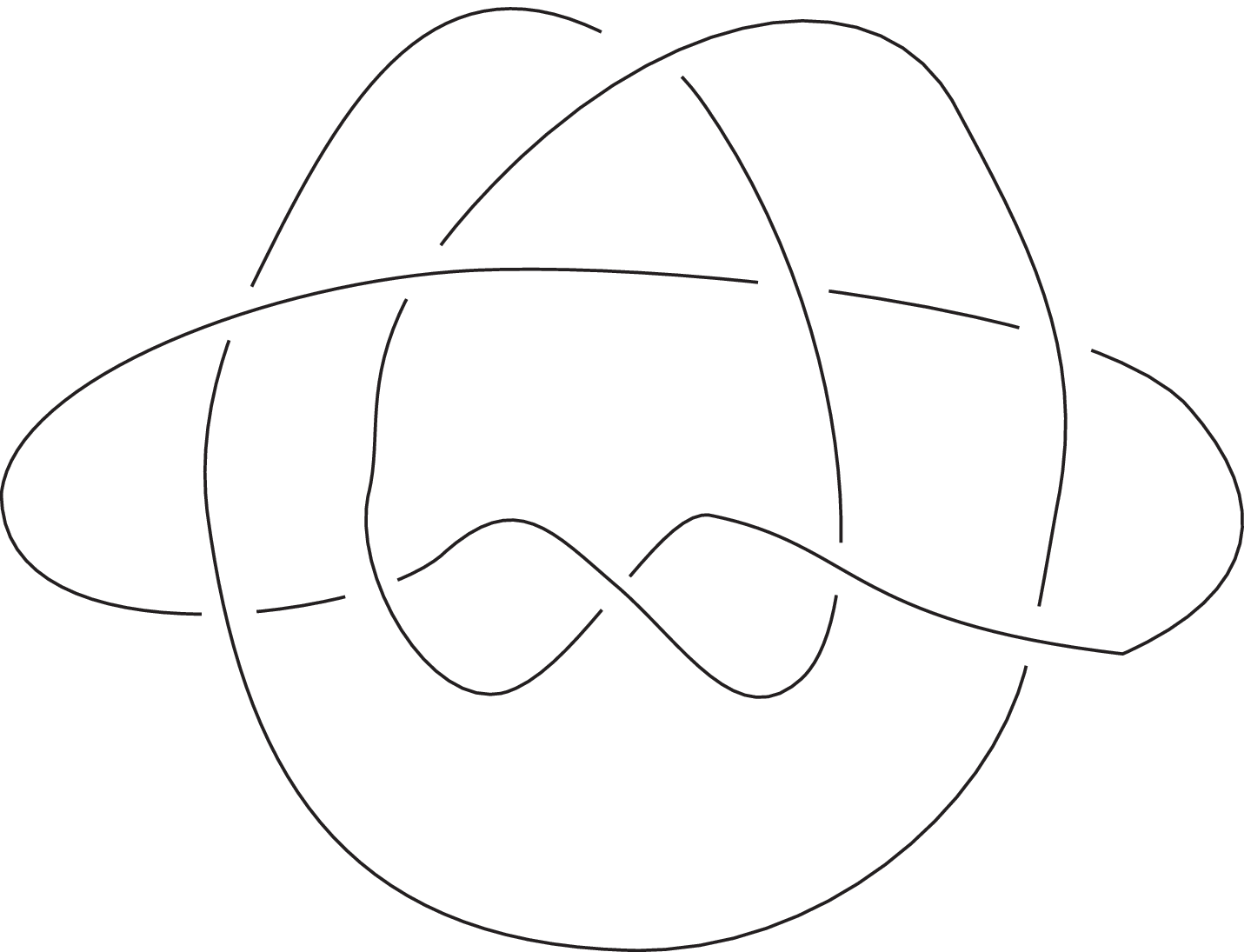}  }
 \caption{${-}10_{161}$} \label{fig10_161}
\end{figure}

\begin{example}\label{10161} The knot ${-}10_{161}$  as illustrated in Figure~\ref{fig10_161} is a 3
stranded braid with 9 positive and 1 negative crossings.  Thus changing one crossing yields a knot with
$\tau = 4$.  This implies that
$\tau({-}10_{161})
\ge 3$.  But, since $g_3({-}10_{161}) = 3$, we also have that
$\tau({-}10_{161})
\le 3$.  So, $\tau({-}10_{161}) = 3$. Since $g_3({-}10_{161}) = 3$ it follows that $g_4(10_{161}) = 3$
also.  (The first calculation of the 4--genus of this knot appeared in~\cite{ta}.)
\end{example}

\begin{figure}[ht!]
\centerline{  \fig{.35}{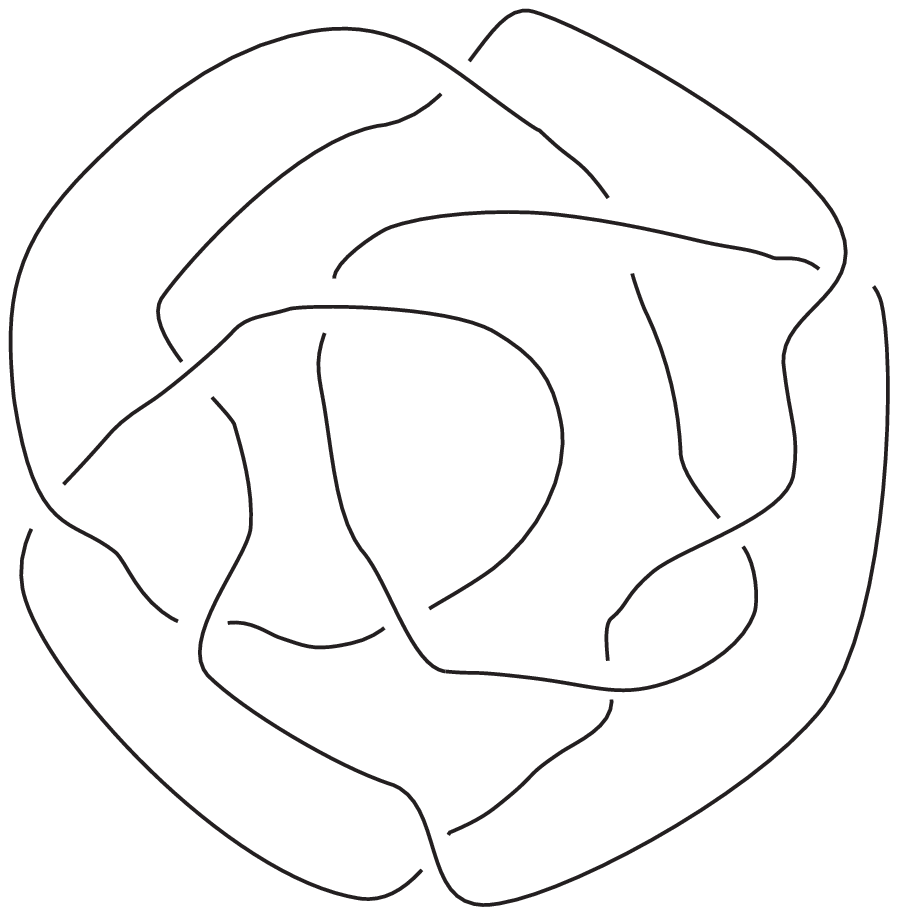}  }
 \caption{$10_{145}$} \label{fig10_145}
\end{figure}

\begin{example}\label{10145} The knot $10_{145}$ is discussed in~\cite{ow}.
 Changing orientation, ${-}10_{145}$ can be drawn as a braid with 4 strands and 11 crossings, 9 of which
are positive.  This is illustrated in Figure~\ref{fig10_145}. Changing two crossing yield a knot
$K$ which by Corollary~\ref{posbraid} has
$\tau(K) = 4$ and so $\tau({-}10_{145}) \ge 2$. On the other hand,
${-}10_{145}$ can be unknotted with 2 crossing changes, so
$\tau({-}10_{145}) \le 2$.  Hence,
$\tau({-}10_{145}) = 2$.  Since the unknotting number of $10_{145}$ is
  at most 2, it follows that $g_4(10_{145}) = 2$. (The first calculation of the 4--genus of this knot also appeared in~\cite{ta}.)

\end{example}

Examples~\ref{10161} and~\ref{10145}  illustrate the following result, extending
Corollary~\ref{posbraid}.   Its proof   follows the exact same lines as the computations in those
examples. 

\begin{corollary} 
  If $\hat{\beta}$ is the closure of a
 braid of
$n$ strands with $k_+$ positive crossings and $k_-$ negative crossings ($k_+ > k_-$), then
$\tau(\hat{\beta}) \ge \frac{k_+ - k_- -n+1}{2}$. 

\end{corollary}

This corollary  immediately gives the bound  $g_4(\hat{\beta}) \ge \frac{k_+ - k_- -n+1}{2}$, the {\it
Slice--Bennequin Inequality} first proved by Rudolph in~\cite{r1}.

\section{Thurston--Bennequin Numbers}  
\begin{figure}[ht!]
\centerline{  \fig{.5}{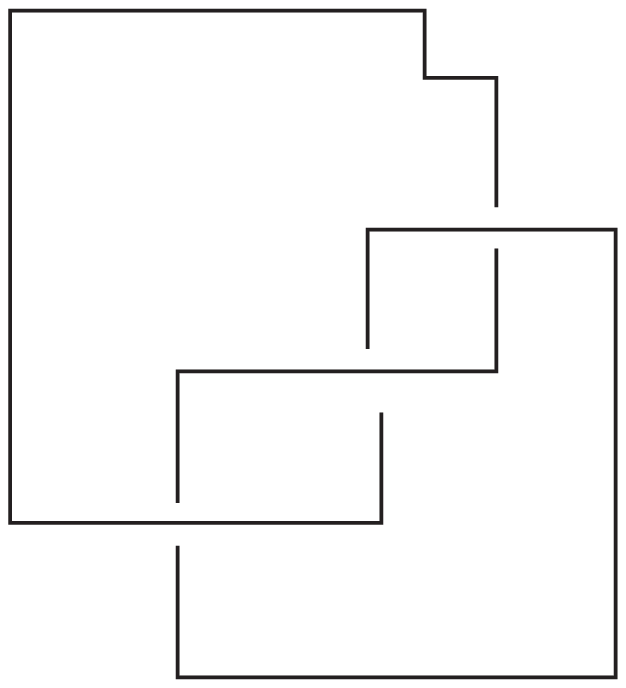}  }
\nocolon \caption{} \label{tref}
\end{figure}

Every knot has a polygonal diagram $D$ consisting of only
vertical and horizontal segments, with each horizontal segment passing over the vertical.  Corners in
such a diagram are naturally labelled northeast, etc.  As described in~\cite{r2}, the
Thurston--Bennequin number of such a diagram, tb($D)$, is the difference of the writhe of the diagram
and the number of northeast corners.  Figure~\ref{tref} illustrates a diagram $D$ of the trefoil knot
with tb($D) = 0$. (In defining  the Thurston--Bennequin number, one usually   considers   knot diagrams
in which all crossings are left handed with respect to the vertical direction and then takes the
difference of the writhe and    the number of right cusps.  The definition here is simply obtained by
``rotating'' the standard definition by 45 degrees.)  The Thurston--Bennequin number of a knot $K$,
TB($K$), is the maximum value of this quantity over all such diagrams for $K$.

\begin{theorem} If the Thurston--Bennequin number of a knot satisfies TB$(K) \ge 0$ then all iterated
untwisted (positive) Whitehead doubles of
$K$,
$Wh_n(K)$,  satisfy
$\tau(Wh_n(K)) = 1$ and thus $g_4(Wh_n(K)) =1$.
\end{theorem}

\begin{proof} Any    polygonal diagram $D$ as above can be isotoped to a diagram
$D'$ with tb($D')= \mbox{tb}(D) -1$; just add a new northeast corner without introducing any new
crossings.  (In~Figure~\ref{tref} an extra northeast corner was   added to a diagram of the trefoil to
illustrate this process.) Thus, we assume tb($D) = 0$.  As observed by Rudolph~\cite{r2}, this diagram
quickly  yields a placement of
$K$ on the fiber $F$ of a torus knot for which the parallel copy $K'$ of $K$ on $F$ has link($K,K') =
0$.    From this   (eg, as in~\cite{r2}) one sees there is also an unknotted curve $\alpha$ on $F$
meeting $K$ transversely in one point, with induced framing
${-}1$.  A neighborhood $G$ of $K \cup
\alpha$ on $F$ is seen to be a genus 1 Seifert surface for the positive untwisted double of $K$, and
hence by Theorem~\ref{mainthm}, $\tau(Wh_1(K)) =1$.  (For this argument to work, one must in fact be a
bit careful in the initial choice of $T_{p,q}$ and $F$; for some choices $K$ embeds, but not
$\alpha$.  Details can be found in the work of Rudolph.)

As shown in~\cite{am} and~\cite{r2}, a simple diagram reveals that  TB($Wh_1(K)) \ge 1$.  Thus, this
process can be iterated.
\end{proof}

We close this section by noting  that in independent work from that presented here,  Olga Plamenevskaya~\cite{pl} has described connections between the Ozsv\'ath--Szab\'o theory, Thurston--Bennequin invariants, and their relationship to the 4--ball genus.

\np

\end{document}
